\documentclass[12pt]{article}
\usepackage{epsfig}
\usepackage{algorithm,algorithmic}
\usepackage{amsmath}
\usepackage{amssymb}
\usepackage{url}
\begin{document}

\def\proof{\noindent {\sl Proof.\ \ }}
\def\endproof{\hfill$\square$
\medskip}

\def\tr{{\rm tr}}
\begin{titlepage}
\title{\vspace{-17mm}
A Purely Algebraic Justification of the
Kabsch-Umeyama Algorithm
%Solution by Singular Value
%Decomposition to the Constrained Orthogonal Procrustes Problem
}
\author{Jim Lawrence$^{1,2}$, Javier Bernal$^1$, Christoph Witzgall$^1$\\
$^1${\small \sl National Institute of Standards and Technology,} \\
{\small \sl Gaithersburg, MD 20899, USA} \\
$^2${\small \sl George Mason University,} \\
{\small \sl 4400 University Dr, Fairfax, VA 22030, USA} \\
{\tt\small $\{$james.lawrence,javier.bernal$\}$@nist.gov \ \ lawrence@gmu.edu}}
%{\small \sl National Institute of Standards and Technology,} \\
%{\small \sl Gaithersburg, MD 20899, USA}}
\date{\ }
\end{titlepage}
\maketitle
\begin{abstract}
%The constrained orthogonal Procrustes problem is the least-squares problem that
%calls for a rotation matrix that optimally aligns two matrices of the same order.
%The algorithm of choice for solving this problem has been the Kabsch-Umeyama algorithm,
%which is effectively no more than the computation of the singular value decomposition
%of a particular matrix. Its justification, as presented separately by Kabsch and Umeyama,
%is not totally algebraic because it is based on solving the minimization problem via Lagrange
%multipliers. In order to provide a more transparent alternative, we present a purely
%algebraic justification of the algorithm.\\[0.2cm]
%
The constrained orthogonal Procrustes problem is the least-squares problem that
calls for a rotation matrix that optimally aligns two matrices of the same order.
Over past decades, the algorithm of choice for solving this
problem has been the Kabsch-Umeyama algorithm, which is effectively no
more than the computation of the singular value decomposition of a particular matrix.
Its justification, as presented separately by Kabsch and Umeyama, is not totally
algebraic since it is based on solving the minimization problem via Lagrange multipliers.
In order to provide a more transparent alternative, it is the main purpose
of this paper to present a purely algebraic justification of the algorithm
through the exclusive use of simple concepts from linear algebra. For the sake
of completeness, a proof is also included of the well known and widely used fact
that the orientation-preserving rigid motion problem,
i.e., the least-squares problem that calls for an orientation-preserving rigid motion
that optimally aligns two corresponding sets of points in $d-$dimensional Euclidean space,
reduces to the constrained orthogonal Procrustes problem.\\[0.2cm]
\end{abstract}
%\par {\small \bf Key words.}
%{\small constrained; Frobenius; least-squares; orientation-preserving; orthogonal; Procrustes;
%rigid motion; rotation; singular value decomposition; trace}
%{\par \small \bf MSC.} {\small 15A18, 15A42, 65H17, 65K99, 93B60}
\section{\large Introduction}
In the {\em orthogonal Procrustes problem} \cite{gower, schonemann}, given real matrices $P$
and $Q$ of size \mbox{$d\times n$}, the problem is that of finding a $d\times d$
orthogonal matrix $U$ that minimizes $\|UQ-P\|_F$, where $\|\cdot\|_F$ denotes the Frobenius
norm of a~matrix. On the other hand, in the {\em constrained orthogonal Procrustes problem}
\cite{kabsch1,kabsch2,umeyama}, the same function is minimized but $U$ is constrained to be
a rotation matrix, i.e., an orthogonal
matrix of determinant~1. By letting $p_i$, $q_i$, $i=1\ldots,n$, be the vectors in
$\mathbb{R}^d$ that are the columns from left to right of $P$ and $Q$, respectively,
since clearly $\|UQ-P\|_F^2=\sum_{i=1}^n \|Uq_i-p_i\|^2$, where \mbox{$\|\cdot\|$} denotes
the $d-$dimensional Euclidean norm, then an alternative formulation of the two problems
above is that of finding an orthogonal matrix $U$ (of determinant~1 for the constrained
problem) that minimizes $\sum_{i=1}^n \|Uq_i-p_i\|^2$. We note that minimizing matrices
do exist for the two problems as the function being minimized is continuous and both the set
of orthogonal matrices and the set of rotation matrices are compact (in some topology).
Finally, in the same vein,
another problem of interest is the {\em orientation-preserving rigid motion problem} which is
that of finding an orientation-preserving rigid motion $\phi$ of~$\mathbb{R}^d$
that minimizes $\sum_{i=1}^n \|\phi(q_i)-p_i\|^2$. An affine linear function $\phi$,
$\phi:\mathbb{R}^d\rightarrow \mathbb{R}^d$, is a rigid motion of $\mathbb{R}^d$ if it is of
the form $\phi(q)=Uq+t$ for $q\in\mathbb{R}^d$, where $U$ is a $d\times d$ orthogonal matrix,
and $t$ is a vector in~$\mathbb{R}^d$. The rigid motion $\phi$ is orientation preserving if
det$(U)=1$, i.e., the determinant of~$U$ equals~1. With $\bar{p}$, $\bar{q}$ denoting the centroids
of $\{p_i\}$, $\{q_i\}$, respectively, as will be shown in Section~3 of this paper, this problem
can be reduced to the constrained orthogonal Procrustes problem by translating $\{p_i\}$,
$\{q_i\}$ to become $\{p_i-\bar{p}\}$, $\{q_i-\bar{q}\}$, respectively, so that the centroid
of each set becomes~$0\in\mathbb{R}^d$.
\par
With $P$, $Q$, $p_i$, $q_i$, $i=1,\ldots,n$, as above, in this paper we focus our attention
mostly on the constrained orthogonal Procrustes problem, and therefore wish to find a
$d\times d$ rotation matrix~$U$ that minimizes $\sum_{i=1}^n \|Uq_i-p_i\|^2$.
\par With this purpose in mind, we rewrite $\sum_{i=1}^n \|Uq_i-p_i\|^2$ as follows, where
given a square matrix $R$, $\tr(R)$ stands for the trace of~$R$.
%\begin{linenomath*}
\begin{eqnarray*}
&& \sum_{i=1}^n ||Uq_i -p_i||^2
= \sum_{i=1}^n (Uq_i - p_i)^T (Uq_i - p_i)
= \tr\bigl((UQ - P)^T(UQ-P)\bigr)\\
\end{eqnarray*}
\begin{eqnarray*}
&=& \tr\bigl((Q^TU^T - P^T) (UQ - P)\bigr)
=\tr(Q^TQ + P^TP - Q^TU^TP - P^T UQ)\\
&=& \tr(Q^TQ) + \tr(P^TP) - 2\tr(P^TUQ).
\end{eqnarray*}
%\end{linenomath*}
Since only the third term in the last line above depends on $U$, it suffices to find
a $d\times d$ rotation matrix $U$ that maximizes $\tr(P^TUQ)$.
Since $\tr(P^TUQ) = \tr(UQP^T)$ (note in general $\tr(AB)=\tr(BA)$, $A$ an $n\times d$ matrix,
$B$ a $d\times n$ matrix), denoting the $d\times d$ matrix $QP^T$ by $M$,
this problem is equivalent to finding a $d\times d$ rotation matrix~$U$ that maximizes
$\tr(UM)$, and it is well known that one such $U$ can be computed from
the singular value decomposition of~$M$~\cite{kabsch1,kabsch2,umeyama}.
This is done with the Kabsch-Umeyama algorithm~\cite{kabsch1,kabsch2,umeyama}
(see Algorithm Kabsch-Umeyama
%Algorithm~\ref{A:kabsch}
below, where $\mathrm{diag}\{s_1,\ldots,s_d\}$
is the $d\times d$ diagonal matrix with numbers $s_1,\ldots,s_d$ as the elements
of the diagonal, in that order running from the upper left to the lower right of the
matrix). A {\em singular value decomposition} (SVD)~\cite{lay} of~$M$
is a representation of the form $M=VSW^T$, where $V$ and $W$ are~$d\times d$ orthogonal
matrices and $S$ is a $d\times d$ diagonal matrix with the singular values of $M$, which
are nonnegative real numbers, appearing in the diagonal of $S$ in descending order,
from  the upper left to the lower right of~$S$. Finally, note that any matrix, not
necessarily square, has a singular value decomposition, not necessarily unique~\cite{lay}.
%\begin{algorithm}
%\caption{Computing optimal $R$}
%\label{A:kabsch}
\begin{algorithmic}
\STATE \noindent\rule{13cm}{0.4pt}
\STATE {\bf Algorithm Kabsch-Umeyama}
\STATE \noindent\rule[.1in]{13cm}{0.4pt}
\STATE Compute $d\times d$ matrix $M=QP^T$.
\STATE Compute SVD of $M$, i.e., identify $d\times d$ matrices $V$, $S$, $W$,
\STATE so that $M = V S W^T$ in the SVD sense.
\STATE Set $s_1= \ldots = s_{d-1}=1$.
\STATE If $\det(VW) > 0$, then set $s_d=1$, else set $s_d=-1$.
\STATE Set $\tilde{S} = \mathrm{diag}\{s_1,\ldots,s_d\}$.
\STATE Return $d\times d$ rotation matrix $U = W \tilde{S} V^T$.
\STATE \noindent\rule{13cm}{0.4pt}
\end{algorithmic}
%\end{algorithm}
\smallskip
\par Algorithm Kabsch-Umeyama has existed for several decades \cite{kabsch1,kabsch2,umeyama},
however the known justifications of the algorithm \cite{kabsch1,kabsch2,umeyama}
are not totally algebraic as they are based on exploiting the optimization technique of
Lagrange multipliers. It is the main purpose of this paper to justify the algorithm in a purely
algebraic manner through the exclusive use of simple concepts from linear algebra.
This is done in Section~2 of the paper. Finally, we note that applications of the algorithm
can be found, notably in the field of functional and shape data
analysis~\cite{dogan,srivastava}, where, in particular, the shapes of two curves are
compared, in part by optimally rotating one curve to match the other.
\section{\large Algebraic Justification of the Kabsch-Umeyama Algorithm}
We justify Algorithm Kabsch-Umeyama using exclusively simple concepts from linear
algebra, mostly in the proof of the following useful proposition.
We note that most of the proof of the proposition is concerned with proving (3) of the
proposition. Thus, it seems reasonable to say that any justification of the algorithm
that requires the conclusion in (3) but lacks a proof for it, is not exactly complete.
See page 47 of the otherwise excellent thesis in~\cite{papa} for an example of this
situation. See~\cite{higham} for an outline of this~dissertation.
\\ \smallskip\\
{\bf Proposition 1:} If $D=\mathrm{diag}\{\sigma_1,\ldots,\sigma_d\}$, $\sigma_j\geq 0$,
$j=1,\ldots,d$, and $W$ is a $d\times d$ orthogonal matrix, then\\
1. tr$(WD)\leq\sum_{j=1}^d \sigma_j$.\\
2. If $B$ is a $d\times d$ orthogonal matrix, $S=B^TDB$, then $\mathrm{tr}(WS)\leq \mathrm{tr}(S)$.\\
3. If det$(W)=-1$, $\sigma_d\leq \sigma_j$, $j=1,\ldots,d-1$, then
$\mathrm{tr}(WD)\leq\sum_{j=1}^{d-1}\sigma_j-\sigma_d$.
\\ \smallskip\\
{\bf Proof:} Since $W$ is orthogonal and if $W_{kj}$, $k,j=1,\ldots,d$, are the
entries of $W$, then, in particular, $W_{jj}\leq 1$, $j=1,\ldots,d$, so that\\
$\mathrm{tr}(WD)=\sum_{j=1}^d W_{jj}\sigma_j\leq \sum_{j=1}^d \sigma_j$, and therefore statement
(1) holds.\\
Accordingly, assuming $B$ is a $d\times d$ orthogonal matrix, since $BWB^T$ is also
orthogonal, it follows from (1) that\\
$\mathrm{tr}(WS)=\mathrm{tr}(WB^TDB)=\mathrm{tr}(BWB^TD)\leq\sum_{j=1}^d\sigma_j=\mathrm{tr}(D)
=\mathrm{tr}(S)$, and therefore (2) holds.\\
If det$(W)=-1$, we show next that a $d\times d$ orthogonal matrix $B$ can be
identified so that with $\bar{W}=B^TWB$, then
$\bar{W}= \left( \begin{smallmatrix}
W_0 & O\\ O^T & -1\\
\end{smallmatrix} \right)$,
$W_0$ interpreted as the upper leftmost $d-1\times d-1$ entries of $\bar{W}$ and as a
$d-1\times d-1$ matrix as well; $O$ interpreted as a vertical column or vector of $d-1$ zeroes.\\
With $I$ as the $d\times d$ identity matrix, then det$(W)=-1$ implies
$\mathrm{det}(W+I)=-\mathrm{det}(W)\mathrm{det}(W+I)=-\mathrm{det}(W^T)\mathrm{det}(W+I)=
-\mathrm{det}(I+W^T)=-\mathrm{det}(I+W)$ which implies det$(W+I)=0$
so that $x\not=0$ exists in $\mathbb{R}^d$ with $Wx=-x$. It also follows then that
$W^TWx=W^T(-x)$ which gives $x=-W^Tx$ so that $W^Tx=-x$ as well.\\
Letting $b_d=x$, vectors $b_1,\ldots,b_{d-1}$ can be obtained so that $b_1,\ldots,b_d$ form
a basis of~$\mathbb{R}^d$, and by the Gram-Schmidt process starting with $b_d$, we may
assume $b_1,\ldots,b_d$ form an orthonormal basis of $\mathbb{R}^d$ with $Wb_d=W^Tb_d=-b_d$.
Letting $B=(b_1,\ldots,b_d)$, interpreted as a $d\times d$ matrix with columns $b_1,\ldots,b_d$,
in that order, it then follows that $B$ is orthogonal, and with $\bar{W}=B^TWB$ and
$W_0$, $O$ as previously described, noting
$B^TWb_d=B^T(-b_d)= \left( \begin{smallmatrix} O\\ -1\\ \end{smallmatrix} \right)$ and
$b_d^TWB=(W^Tb_d)^TB=(-b_d)^TB=(O^T \, -1)$,
%left( \begin{smallmatrix} 0^T & -1\\ \end{smallmatrix} \right)$
then $\bar{W}= \left( \begin{smallmatrix}
W_0 & O\\ O^T & -1\\
\end{smallmatrix} \right)$.
Note $\bar{W}$ is orthogonal and therefore so is the $d-1\times d-1$ matrix~$W_0$.\\
Let $S=B^TDB$ and write
$S = \left( \begin{smallmatrix}
S_0 & a\\ b^T & \gamma\\
\end{smallmatrix} \right)$,
$S_0$ interpreted as the upper leftmost $d-1\times d-1$ entries of $S$ and as a
$d-1\times d-1$ matrix as well; $a$ and $b$ interpreted as vertical columns or vectors of
$d-1$ entries, and $\gamma$ as a scalar.\\
Note $\mathrm{tr}(WD)=\mathrm{tr}(B^TWDB)=\mathrm{tr}(B^TWBB^TDB)=\mathrm{tr}(\bar{W}S)$,
so that $\bar{W}S=$
$\left( \begin{smallmatrix}
W_0 & O\\ O^T & -1\\
\end{smallmatrix} \right)$
$\left( \begin{smallmatrix}
S_0 & a\\ b^T & \gamma\\
\end{smallmatrix} \right)=$
$\left( \begin{smallmatrix}
W_0S_0 & W_0a\\ -b^T & -\gamma\\
\end{smallmatrix} \right)$
gives $\mathrm{tr}(WD)=\mathrm{tr}(W_0S_0)-\gamma$.\\
We show $\mathrm{tr}(W_0S_0)\leq\mathrm{tr}(S_0)$. For this purpose let
$\hat{W}= \left( \begin{smallmatrix}
W_0 & O\\ O^T & 1\\
\end{smallmatrix} \right)$,
$W_0$ and $O$ as above. Since $W_0$ is orthogonal, then clearly $\hat{W}$ is a $d\times d$
orthogonal matrix, and by (2), $\mathrm{tr}(\hat{W}S)\leq \mathrm{tr}(S)$
so that $\hat{W}S=$
$\left( \begin{smallmatrix}
W_0 & O\\ O^T & 1\\
\end{smallmatrix} \right)$
$\left( \begin{smallmatrix}
S_0 & a\\ b^T & \gamma\\
\end{smallmatrix} \right)=$
$\left( \begin{smallmatrix}
W_0S_0 & W_0a\\ b^T & \gamma\\
\end{smallmatrix} \right)$
gives $\mathrm{tr}(W_0S_0)+\gamma=\mathrm{tr}(\hat{W}S)\leq\mathrm{tr}(S)=\mathrm{tr}(S_0)+\gamma$.
Thus, $\mathrm{tr}(W_0S_0)\leq\mathrm{tr}(S_0)$.\\
Note $\mathrm{tr}(S_0)+\gamma = \mathrm{tr}(S)= \mathrm{tr}(D)$, and if $B_{kj}$,
$k,j=1,\ldots,d$ are the entries of~$B$, then $\gamma = \sum_{k=1}^d B_{kd}^2\sigma_k$,
a convex combination of the $\sigma_k$'s, so that $\gamma\geq\sigma_d$.
It then follows that\\ $\mathrm{tr}(WD)=\mathrm{tr}(W_0S_0)-\gamma\leq\mathrm{tr}(S_0)-\gamma=
\mathrm{tr}(D)-\gamma-\gamma\leq\sum_{j=1}^{d-1}\sigma_j-\sigma_d$, and therefore
(3) holds.
\endproof
%\ $\Box$
\\ \smallskip
\par
Finally, the following theorem, a consequence of Proposition 1, justifies the Kabsch-Umeyama
algorithm.
\\ \smallskip\\
{\bf Theorem 1:} Given a $d\times d$ matrix $M$, let $V$, $S$, $W$ be $d\times d$ matrices such that
the singular value decomposition of $M$ gives $M=VSW^T$. If det$(VW)>0$, then $U=WV^T$
maximizes~tr$(UM)$ over all $d\times d$ rotation matrices~$U$.
Otherwise, if det$(VW)<0$, with $\tilde{S}=\mathrm{diag}\{s_1,\ldots,s_d\}$,
$s_1 = \ldots = s_{d-1}=1,\ s_d=-1$, then $U=W\tilde{S}V^T$ maximizes~tr$(UM)$ over all $d\times d$
rotation matrices~$U$.
\\ \smallskip\\
{\bf Proof:} Let $\sigma_j$, $j=1,\ldots,d$, $\sigma_1\geq \sigma_2\geq\ldots\geq\sigma_d\geq 0$,
be the singular values of~$M$, so that $S=\mathrm{diag}\{\sigma_1,\ldots,\sigma_d\}$.\\
Assume det$(VW)>0$. If $U$ is any rotation matrix, then $U$ is orthogonal.
From (1) of Proposition~1 since $W^TUV$ is orthogonal, then\\
$\mathrm{tr}(UM)=\mathrm{tr}(UVSW^T)=\mathrm{tr}(W^TUVS)\leq\sum_{j=1}^d\sigma_j.$\\
On the other hand, if $U=WV^T$, then $U$ is clearly orthogonal, det$(U)=1$, and
$\mathrm{tr}(UM)=\mathrm{tr}(WV^TVSW^T)=\mathrm{tr}(WSW^T)=\mathrm{tr}(S)
=\sum_{j=1}^d\sigma_j.$\\
Thus, $U=WV^T$ maximizes~tr$(UM)$ over all $d\times d$ rotation matrices~$U$.\\
Finally, assume det$(VW)<0$. If $U$ is any rotation matrix, then $U$ is orthogonal
and~det$(U)=1$. From (3) of Proposition~1 since $W^TUV$ is orthogonal and det$(W^TUV)=-1$, then\\
$\mathrm{tr}(UM)=\mathrm{tr}(UVSW^T)=\mathrm{tr}(W^TUVS)\leq\sum_{j=1}^{d-1}\sigma_j-\sigma_d.$\\
On the other hand, if $U=W\tilde{S}V^T$, then $U$ is clearly orthogonal, det$(U)=1$, and
$\mathrm{tr}(UM)=\mathrm{tr}(W\tilde{S}V^TVSW^T)=\mathrm{tr}(W\tilde{S}SW^T)=
\mathrm{tr}(\tilde{S}S) =\sum_{j=1}^{d-1}\sigma_j-\sigma_d.$\\
Thus, $U=W\tilde{S}V^T$ maximizes~tr$(UM)$ over all $d\times d$ rotation matrices~$U$.
\endproof
%\ $\Box$
\section{\large Reduction of the Orientation-Preserving Rigid Motion Problem to the
Constrained Orthogonal Procrustes Problem}
Although not exactly related to the main goal of this paper,
for the sake of completeness, we show the orientation-preserving rigid motion problem
reduces to the constrained orthogonal Procrustes problem.
For this purpose, let $\bar q$ and $\bar p$ denote the centroids of the sets $\{q_i\}_{i=1}^n$
and $\{p_i\}_{i=1}^n$ in~$\mathbb{R}^d$, respectively:
%\begin{linenomath*}
$$\bar q = \frac{1}{n} \sum_{i=1}^n q_i\ \ \ \mathrm{and}
\ \ \ \bar p = \frac{1}{n} \sum_{i=1}^n p_i\ . $$
%\end{linenomath*}
First we prove a proposition that shows, in particular, that if $\hat{\phi}(\bar{q})\not = \bar{p}$, 
then $\phi=\hat{\phi}$ does not minimize
%\begin{linenomath*}
$$\Delta(\phi)=\sum_{i=1}^n \|\phi(q_i)-p_i\|^2,$$
%\end{linenomath*}
the minimization occurring
over either the set
of all rigid motions $\phi$ of~$\mathbb{R}^d$ or the smaller set of rigid motions $\phi$
of~$\mathbb{R}^d$ that are orientation preserving.
\\ \smallskip\\
{\bf Proposition 2:} Let $\phi$ be a rigid motion of $\mathbb{R}^d$ with
$\phi(\bar{q})\not=\bar{p}$ and define an affine linear function $\tau$,
$\tau: \mathbb{R}^d\rightarrow\mathbb{R}^d$, by $\tau(q) = \phi(q)-\phi(\bar{q})+\bar{p}$
for $q\in\mathbb{R}^d$. Then $\tau$ is a rigid motion of $\mathbb{R}^d$, $\tau(\bar{q})=\bar{p}$,
$\Delta(\tau) < \Delta(\phi)$, and if $\phi$ is
orientation preserving, then so is~$\tau$.\\ \smallskip\\
{\bf Proof:} Clearly $\tau(\bar{q})=\bar{p}$. Let $U$ be a $d\times d$ orthogonal matrix
and $t \in \mathbb{R}^d$ such that $\phi(q) = Uq + t$ for $q\in\mathbb{R}^d$.
Then $\tau(q) = Uq - U\bar q + \bar p$ so that $\tau$ is a rigid motion of~$\mathbb{R}^d$,
$\tau$ is orientation preserving if $\phi$ is, and for $ 1 \leq i \leq n$, we~have
\begin{eqnarray*}
&& ||\phi(q_i) - p_i||^2 - ||\tau(q_i) - p_i||^2
= (Uq_i + t - p_i)^T(Uq_i + t - p_i)\\
%\end{eqnarray*}
%\begin{eqnarray*}
&& - (Uq_i - U\bar q + \bar p - p_i)^T (Uq_i - U\bar q + \bar p - p_i) \\
&=& \big((Uq_i - p_i)^T(Uq_i - p_i) + 2(Uq_i - p_i)^T t + t^T t \big) -
\big((Uq_i - p_i)^T(Uq_i - p_i)\\ && - 2(Uq_i - p_i)^T(U\bar q - \bar p)
 + (U \bar q - \bar p)^T(U \bar q - \bar p)\big)\\
&=& 2(Uq_i - p_i + t)^T(U\bar q - \bar p + t) - (U \bar q - \bar p + t)^T
(U \bar q - \bar p + t).
\end{eqnarray*}
%It then follows that $\Delta(\phi) - \Delta(\tau)$ equals
It then follows that
\begin{eqnarray*}
&& \Delta(\phi) - \Delta(\tau)\\
&=& \sum_{i=1}^n \big(2(Uq_i - p_i + t)^T(U\bar q - \bar p + t) - (U \bar q - \bar p + t)^T
(U \bar q - \bar p + t)\big)\\
&=& n||U\bar q - \bar p + t||^2=n||\phi(\bar q) - \bar p||^2 >0
%\mathrm{\ as\ } \phi(\bar q) - \bar p \mathrm{\ is\ nonzero.}
%\hspace{.6in}  \Box
\end{eqnarray*}
as $\phi(\bar q) - \bar p$  is nonzero. Thus $\Delta(\tau)<\Delta(\phi)$.
\endproof
\\ \smallskip
\par Finally, the following corollary, a consequence of Proposition 2, shows that the
problem of finding an orientation-preserving rigid motion $\phi$ of~$\mathbb{R}^d$ that minimizes
$\sum_{i=1}^n \|\phi(q_i)-p_i\|^2$ can be reduced to a constrained orthogonal Procrustes
problem which, of course, then can be solved with the Kabsch-Umeyama algorithm. Here
$r_i=p_i-\bar{p}$, $s_i=q_i-\bar{q}$, for $i=1,\ldots,n$, and if
$\bar{r}=\frac{1}{n}\sum_{i=1}^nr_i$, $\bar{s}=\frac{1}{n}\sum_{i=1}^ns_i$,
then clearly $\bar{r}=\bar{s}=0$.
\\ \smallskip\\
{\bf Corollary 1:} Let $\hat{U}$ be such that $U=\hat{U}$ minimizes
$\sum_{i=1}^n \|Us_i-r_i\|^2$ over all $d\times d$ rotation matrices~$U$.
Let $\hat{t}=\bar{p}-\hat{U}\bar{q}$, and let $\hat{\phi}$ be given by
$\hat{\phi}(q)=\hat{U}q+\hat{t}$ for $q\in\mathbb{R}^d$. Then $\phi=\hat{\phi}$ minimizes
$\sum_{i=1}^n \|\phi(q_i)-p_i\|^2$ over all orientation-preserving rigid motions~$\phi$
of~$\mathbb{R}^d$.
\\ \smallskip\\
{\bf Proof:} One such $\hat{U}$ can be computed with the Kabsch-Umeyama algorithm.\\
By Proposition~2, in order to minimize $\sum_{i=1}^n \|\phi(q_i)-p_i\|^2$ over all
orientation-preserving rigid motions $\phi$ of $\mathbb{R}^d$, it suffices to do it over those
for which $\phi(\bar{q})=\bar{p}$. Therefore, it suffices to minimize
$\sum_{i=1}^n \|Uq_i+t-p_i\|^2$ with $t=\bar{p}-U\bar{q}$ over all $d\times d$ rotation
matrices~$U$, i.e., it suffices to minimize
%\begin{linenomath*}
$$\sum_{i=1}^n\|Uq_i+\bar{p}-U\bar{q}-p_i\|^2 =
\sum_{i=1}^n\|(U(q_i-\bar{q})-(p_i-\bar{p})\|^2$$
%\end{linenomath*}
over all $d\times d$ rotation matrices~$U$. But minimizing the last expression is equivalent
to minimizing $\sum_{i=1}^n \|Us_i-r_i\|^2$ over all $d\times d$ rotation matrices~$U$.
Since $U=\hat{U}$ is a solution to this last problem, it then follows that
$U=\hat{U}$ minimizes $\sum_{i=1}^n\|Uq_i+\bar{p}-U\bar{q}-p_i\|^2$
$=\sum_{i=1}^n \|Uq_i+t-p_i\|^2$ with $t=\bar{p}-U\bar{q}$ over all $d\times d$
rotation matrices~$U$. Consequently, 
if $\hat{t}=\bar{p}-\hat{U}\bar{q}$, and $\hat{\phi}$ is given by
$\hat{\phi}(q)=\hat{U}q+\hat{t}$ for $q\in\mathbb{R}^d$,
then $\phi=\hat{\phi}$ clearly minimizes $\sum_{i=1}^n \|\phi(q_i)-p_i\|^2$
over all orientation-preserving rigid motions~$\phi$ of~$\mathbb{R}^d$.
\endproof
%$\ \Box$
%\\ \smallskip\\
%\\ \pagebreak\\
%{\bf\large Summary}\\ \smallskip\\
%In

\end{document}